\documentclass{amsart}
\usepackage{graphicx,color}
\newtheorem{Thm}{Theorem}

\begin{document}

\author{Yeonhee Jang}
\address{Department of Mathematics, Nara Women's University, Kitauoyanishi-machi, Nara 630-8506, Japan}
\email{yeonheejang@cc.nara-wu.ac.jp}

\author{Tsuyoshi Kobayashi}
\address{Department of Mathematics, Nara Women's University, Kitauoyanishi-machi, Nara 630-8506, Japan}
\email{tsuyoshi@cc.nara-wu.ac.jp}

\author{Makoto Ozawa}
\address{Department of Natural Sciences, Faculty of Arts and Sciences, Komazawa University, 1-23-1 Komazawa, Setagaya-ku, Tokyo, 154-8525, Japan}
\email{w3c@komazawa-u.ac.jp}

\author{Kazuto Takao}
\address{Institute of Mathematics for Industry, Kyushu University, 744 Motooka, Nishi-ku, Fukuoka 819-0395, Japan}
\email{takao@imi.kyushu-u.ac.jp}

\title[A knot with destabilized bridge spheres]{A knot with destabilized bridge spheres of arbitrarily high bridge number}
\subjclass[2010]{57M25}
\keywords{knots, bridge spheres, stabilizations}
\thanks{The first author is supported by Grant-in-Aid for Research Activity Start-up (No.\ 25887039).
The second author is supported by Grant-in-Aid for Scientific Research (C) (No.\ 25400091).
The third author is partially supported by Grant-in-Aid for Scientific Research (C) (No.\ 23540105) and Grant-in-Aid for Scientific Research (C) (No.\ 26400097).
The fourth author is supported by Grant-in-Aid for JSPS Fellows (No.\ 26$\cdot $3358) and Grant-in-Aid for Young Scientists (B) (No.\ 26800042)}

\begin{abstract}
We show that there exists an infinite family of knots, each of which has, for each integer $k\geq 0$, a destabilized $(2k+5)$-bridge sphere.
We also show that, for each integer $n\geq 4$, there exists a knot with a destabilized $3$-bridge sphere and a destabilized $n$-bridge sphere.
\end{abstract}

\maketitle

\section{Introduction}\label{introduction}

We study bridge spheres of knots in the $3$-sphere, paying particular attention to destabilized ones.
Any knot admits infinitely many bridge spheres, and to classify them is a general problem.
In particular, studying destabilized ones is essential because all the others can be obtained from them by stabilizations up to isotopy.
(See Section~\ref{preliminaries} for definitions.)

Some knots admit unique destabilized bridge spheres up to isotopy.
Otal \cite{Otal1} showed that the trivial knot admits no other destabilized bridge spheres than those isotopic to the canonical $1$-bridge sphere.
Schubert \cite{Schubert} showed that each non-trivial rational knot admits a unique $2$-bridge sphere up to isotopy, and Otal \cite{Otal2} showed that it admits no other destabilized bridge spheres than those isotopic to the $2$-bridge sphere.
The third author \cite{Ozawa} obtained the corresponding result for torus knots.
Zupan \cite{Zupan} showed that if a meridionally small knot admits a unique destabilized bridge sphere up to isotopy, then any cable knot of it also does.

Some knots, however, admit multiple non-isotopic destabilized bridge spheres.
Birman \cite{Birman2} gave a composite knot with two non-isotopic destabilized $3$-bridge spheres.
Montesinos \cite{Montesinos} gave a prime knot with two non-isotopic destabilized $3$-bridge spheres.
Johnson--Tomova \cite{Johnson-Tomova} gave a knot with two non-isotopic destabilized bridge spheres which are far apart in the sense of stable equivalence.
The first author \cite{Jang1,Jang2} gave a knot with four destabilized $3$-bridge spheres which are pairwise non-isotopic.
The third and fourth authors \cite{Ozawa-Takao} gave a knot with a destabilized $3$-bridge sphere and a destabilized $4$-bridge sphere.

\begin{figure}[ht]
\begin{picture}(275,200)(0,0)
\put(55,0){\includegraphics{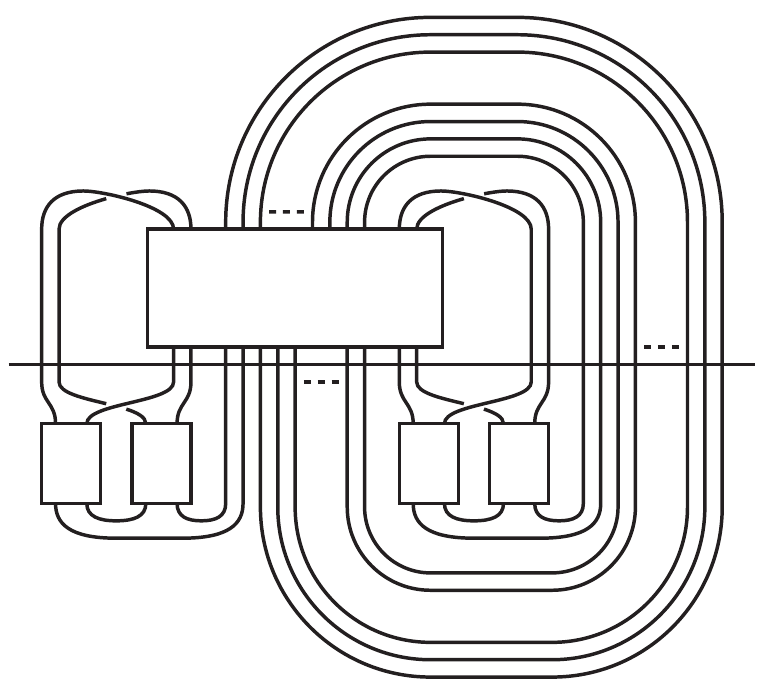}}
\put(117,143){$\overbrace{\hspace{49pt}}^{\colorbox{white}{2k+1}}$}
\put(115,114){$q$ full twists}
\put(72,65){$p_1$}
\put(98,65){$p_2$}
\put(175,65){$p_3$}
\put(201,65){$p_4$}
\put(0,91){$S_{p_1,p_2,p_3,p_4,q,k}$}
\end{picture}
\caption{
A knot $K_{p_1,p_2,p_3,p_4,q,k}$ and a sphere $S_{p_1,p_2,p_3,p_4,q,k}$.
The box labelled ``$p_i$'' represents the $p_i$ left-handed half-twists.
See also the top of Figure~\ref{fig_deformation} for the case where $p_1=p_2=p_3=p_4=2$, $q=0$ and $k=4$.}
\label{fig_quotient}
\end{figure}

In this paper, we give an infinite family of knots, each of which has a $5$-bridge sphere and destabilized bridge spheres of arbitrarily high bridge number, as follows.
\begin{Thm}\label{quotient}
Let $K_{p_1,p_2,p_3,p_4,q,k}$ and $S_{p_1,p_2,p_3,p_4,q,k}$ be the knot and the sphere, respectively, shown in Figure~\ref{fig_quotient}, for integers $p_1,p_2,p_3,p_4,q$ and a non-negative integer $k$.
The knot type of $K_{p_1,p_2,p_3,p_4,q,k}$ does not depend on $k$, and if $|2p_i+1|\geq 5$ for each $i\in \{ 1,2,3,4\} $ and $q$ is an even number with $|q|\geq 12$, then $S_{p_1,p_2,p_3,p_4,q,k}$ is a destabilized $(2k+5)$-bridge sphere of $K_{p_1,p_2,p_3,p_4,q,k}$.
\end{Thm}
\noindent
For the proof of the assertion that $S_{p_1,p_2,p_3,p_4,q,k}$ is destabilized, we show the strong irreducibility of it.

\begin{figure}[ht]
\begin{picture}(267,120)(0,0)
\put(7,10){\includegraphics{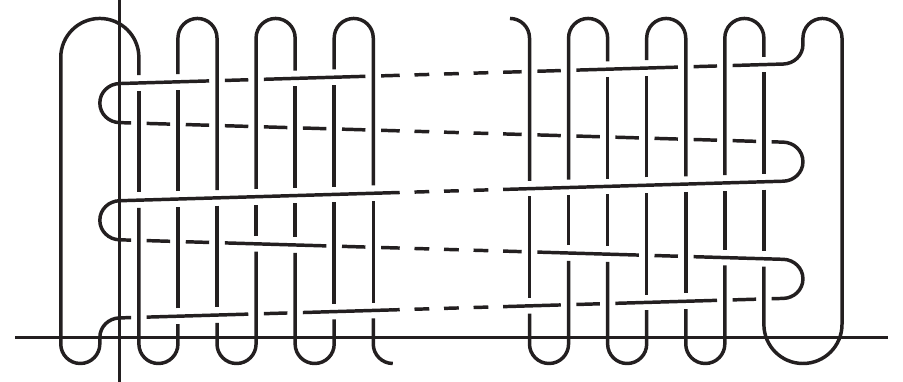}}
\put(0,20){$S_n$}
\put(37,0){$S_n^\perp $}
\end{picture}\\
in the case where $n$ is even\\[12pt]
\begin{picture}(267,120)(0,0)
\put(7,10){\includegraphics{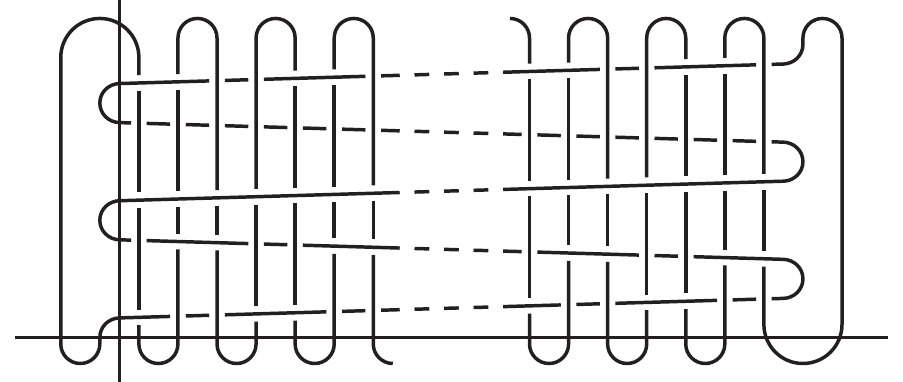}}
\put(0,20){$S_n$}
\put(37,0){$S_n^\perp $}
\end{picture}\\
in the case where $n$ is odd
\caption{
A knot $K_n$ and spheres $S_n^\perp $ and $S_n$.
The knot $K_n$ intersects $S_n$ in $2n$ points.
The second top weft runs under the warps, and the third weft runs over the warps.
Each of the other wefts threads across the warps by repeating ``over, over, under, under'' with the exception in the leftmost and rightmost parts.}
\label{fig_biplat}
\end{figure}

It remains an open problem whether there exists a knot that has a $3$-bridge or $4$-bridge sphere and destabilized bridge spheres of arbitrarily high bridge number.
We remark that such a knot cannot have a $2$-bridge sphere, by the fact that only rational knots may admit $2$-bridge spheres, and the results about rational knots mentioned above.
As an answer for a subordinate problem, we give the following theorem.
\begin{Thm}\label{biplat}
Let $K_n$, $S_n^\perp $ and $S_n$ be the knot and the two spheres, respectively, shown in Figure~\ref{fig_biplat}, for an integer $n$ with $n\geq 4$.
Then, $S_n^\perp $ is a destabilized $3$-bridge sphere of $K_n$, and $S_n$ is a destabilized $n$-bridge sphere of $K_n$.
\end{Thm}
\noindent
For the proof of the assertion that $S_n$ is destabilized, we show the strong irreducibility of it.

This work is motivated by Casson--Gordon's work for Heegaard surfaces of $3$-manifolds.
They gave an infinite family of $3$-manifolds, each of which has destabilized Heegaard surfaces of arbitrarily high genus, and we use their result for the proof of Theorem~\ref{quotient} (see Section~\ref{proof_quotient}).
Also, they introduced a criterion for recognizing destabilized Heegaard surfaces and used it for their proof, and similarly we introduce a criterion for recognizing destabilized bridge spheres (see Section~\ref{diagrams}) and use it for our proof of Theorem~\ref{biplat} (see Section~\ref{proof_biplat}).

\section{Heegaard surfaces and bridge spheres}\label{preliminaries}

In this section, we briefly review standard definitions and facts concerning Heegaard surfaces of $3$-manifolds and bridge spheres of links.

\subsection{Heegaard surfaces}

A {\it handlebody} is a $3$-manifold $H$ homeomorphic to a closed regular neighborhood of a connected finite graph embedded in ${\mathbb R}^3$.
An {\it essential disk} of $H$ is a properly embedded disk $D$ in $H$ such that $\partial D$ does not bound a disk in the surface $\partial H$.

A {\it Heegaard surface} of a closed orientable $3$-manifold $M$ is a closed surface $\Sigma \subset M$ which decomposes $M$ into two handlebodies.
That is to say, there are two handlebodies $H^+,H^-\subset M$ such that $H^+\cup H^-=M$ and $H^+\cap H^-=\partial H^+=\partial H^-=\Sigma $.
Two Heegaard surfaces $\Sigma _1,\Sigma _2$ of $M$ are said to be {\it isotopic} if there is an ambient isotopy $\{ F_t:M\rightarrow M\} _{t\in [0,1]}$ such that $F_0=id_M$ and $F_1(\Sigma _1)=\Sigma _2$.

\begin{figure}[ht]
\begin{picture}(270,80)(0,0)
\put(0,0){\includegraphics{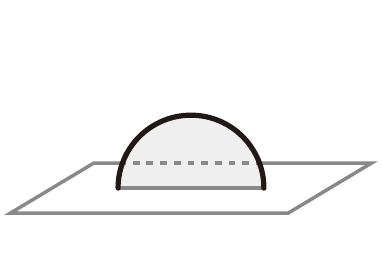}}
\put(125,25){{\Large $\rightarrow $}}
\put(150,0){\includegraphics{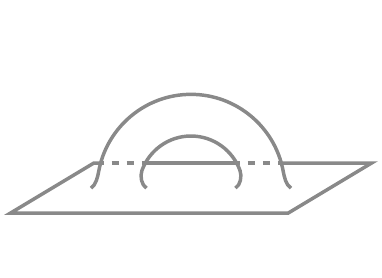}}
\put(0,25){$\Sigma $}
\put(52,53){$\gamma $}
\put(260,25){$\widetilde{\Sigma }$}
\end{picture}
\caption{A stabilization for a Heegaard surface.}
\label{fig_stab_h}
\end{figure}

The notion of stabilization for Heegaard surfaces is defined as follows.
Let $M$, $\Sigma $, $H^+,H^-$ be as above.
Let $\gamma $ be a properly embedded arc in $H^+$ parallel to $\Sigma $, and $N(\gamma )$ be a closed regular neighborhood of $\gamma $.
Let $\widetilde{H}^-$ denote the union $H^-\cup N(\gamma )$, $\widetilde{H}^+$ denote the closure of $H^+\setminus N(\gamma )$, and $\widetilde{\Sigma }$ denote their common boundary (see Figure~\ref{fig_stab_h}).
Then, both $\widetilde{H}^+,\widetilde{H}^-$ are handlebodies and so $\widetilde{\Sigma }$ is a Heegaard surface of $M$.
We say that $\widetilde{\Sigma }$ is obtained from $\Sigma $ by a {\it stabilization}.
Note that the genus of $\widetilde{\Sigma }$ is greater than that of $\Sigma $ by one, and that the isotopy class of $\widetilde{\Sigma }$ depends only on the isotopy class of $\Sigma $.

A Heegaard surface is said to be {\it destabilized} if it cannot be obtained from another Heegaard surface by a stabilization.

A Heegaard surface is said to be {\it strongly irreducible} if any essential disk of one of the handlebodies and any essential disk of the other have non-empty intersection.
It is known that a strongly irreducible Heegaard surface is destabilized, with the exception of Heegaard surfaces of genus $1$ of $S^3$.

\subsection{Bridge spheres}

An {\it $n$-string trivial tangle} is the pair $(B,\tau )$ of a $3$-ball $B$ and a collection $\tau $ of pairwise disjoint properly embedded arcs $\tau _1,\tau _2,\ldots ,\tau _n$ in $B$ simultaneously parallel to $\partial B$.
There are pairwise disjoint disks $E_1,E_2,\ldots ,E_n$ in $B$ such that $E_i$ is cobounded by $\tau _i$ and an arc in $\partial B$ for each $i\in \{ 1,2,\ldots ,n\} $.
We call $\{ E_1,E_2,\ldots ,E_n\} $ a {\it complete collection of bridge disks} of $(B,\tau )$.
An {\it essential disk} of $(B,\tau )$ is a properly embedded disk $D$ in $B$ such that $D\cap \tau =\emptyset $ and $\partial D$ does not bound a disk in the $2n$-punctured sphere $\partial B\setminus \tau $.

An {\it $n$-bridge sphere} of a link $L$ in the $3$-sphere $S^3$ is a $2$-sphere $S\subset S^3$ which is  transverse to $L$ and decomposes $(S^3,L)$ into two $n$-string trivial tangles.
That is to say, the pairs $(B^+,\tau ^+),(B^-,\tau ^-)$ are $n$-string trivial tangles, where $B^+,B^-\subset S^3$ denote the $3$-balls divided by $S$ and $\tau ^+=L\cap B^+,\tau ^-=L\cap B^-$.
We call the number $n$ the {\it bridge number} of the bridge sphere $S$.
Two bridge spheres $S_1,S_2$ of $L$ are said to be {\it isotopic} if there is an ambient isotopy $\{ F_t:S^3\rightarrow S^3\} _{t\in [0,1]}$ such that $F_0=id_{S^3}$, $F_1(S_1)=S_2$, and $F_t(S_1)$ is transverse to $L$ for $t\in[0,1]$.

\begin{figure}[ht]
\begin{picture}(270,80)(0,0)
\put(0,0){\includegraphics{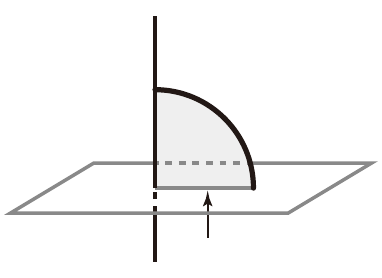}}
\put(125,25){{\Large $\rightarrow $}}
\put(150,0){\includegraphics{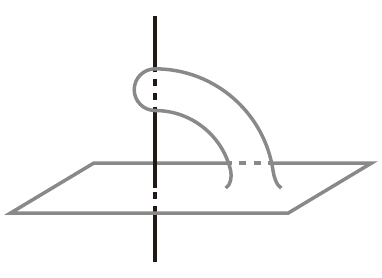}}
\put(32,66){$K$}
\put(0,25){$S$}
\put(35,37){$\alpha $}
\put(56,2){$\beta $}
\put(67,50){$\gamma $}
\put(260,25){$\widetilde{S}$}
\end{picture}
\caption{A stabilization for a bridge sphere.}
\label{fig_stab_b}
\end{figure}

The notion of stabilization for bridge spheres is defined as follows.
Let $L$, $S$, $n$, $B^+$, $B^-$, $\tau ^+$, $\tau ^-$ be as above.
Let $\alpha ,\beta ,\gamma $ be arcs in $\tau ^+,S,B^+$, respectively, such that $\gamma \cap S$ consists of one endpoint of $\gamma $ and $\gamma \cap \tau ^+$ consists of the other endpoint, and $\alpha \cup \beta \cup \gamma $ is the boundary of a disk in $B^+$ whose interior is disjoint from $\tau ^+$.
Let $N(\gamma )$ be a closed regular neighborhood of $\gamma $.
Let $\widetilde{B}^-$ denote the union $B^-\cup N(\gamma )$, $\widetilde{B}^+$ denote the closure of $B^+\setminus N(\gamma )$, and $\widetilde{S}$ denote their common boundary (see Figure~\ref{fig_stab_b}).
Let $\widetilde{\tau }^+=L\cap \widetilde{B}^+$ and $\widetilde{\tau }^-=L\cap \widetilde{B}^-$.
Then, both $(\widetilde{B}^+,\widetilde{\tau }^+),(\widetilde{B}^-,\widetilde{\tau }^-)$ are $(n+1)$-string trivial tangles and so $\widetilde{S}$ is an $(n+1)$-bridge sphere of $L$.
We say that $\widetilde{S}$ is obtained from $S$ by a {\it stabilization}.
Note that the isotopy class of $\widetilde{S}$ depends only on the isotopy class of $S$ if $L$ is a knot.

A bridge sphere is said to be {\it destabilized} if it cannot be obtained from another bridge sphere by a stabilization.

A bridge sphere is said to be {\it strongly irreducible} if any essential disk of one of the trivial tangles and any essential disk of the other have non-empty intersection.
It is known that a strongly irreducible bridge sphere is destabilized, with the exception of $2$-bridge spheres of the trivial knot.

\section{Proof of Theorem~1}\label{proof_quotient}

In this section, we give a proof of Theorem~\ref{quotient}.
To do this, we use the result of Casson--Gordon's about Heegaard surfaces of $3$-manifolds as mentioned in Section~\ref{introduction}.

We use, for simplicity, the symbols $K_{q,k}$ and $S_{q,k}$ for the knot $K_{p_1,p_2,p_3,p_4,q,k}$ and the sphere $S_{p_1,p_2,p_3,p_4,q,k}$, respectively, shown in Figure~\ref{fig_quotient}.
At this stage, we make no assumption on the integers $p_1,p_2,p_3,p_4,q$ and the non-negative integer $k$.

First of all, it is easy to see that the sphere $S_{q,k}$ is a $(2k+5)$-bridge sphere of the knot $K_{q,k}$.
Indeed, $K_{q,k}$ is in a $(2k+5)$-bridge position with respect to the height in Figure~\ref{fig_quotient}, and $S_{q,k}$ is a dividing sphere for $K_{q,k}$ (see \cite{Scharlemann} for the concept of bridge position).

\begin{figure}[ht]
\begin{picture}(220,200)(0,0)
\put(0,0){\includegraphics{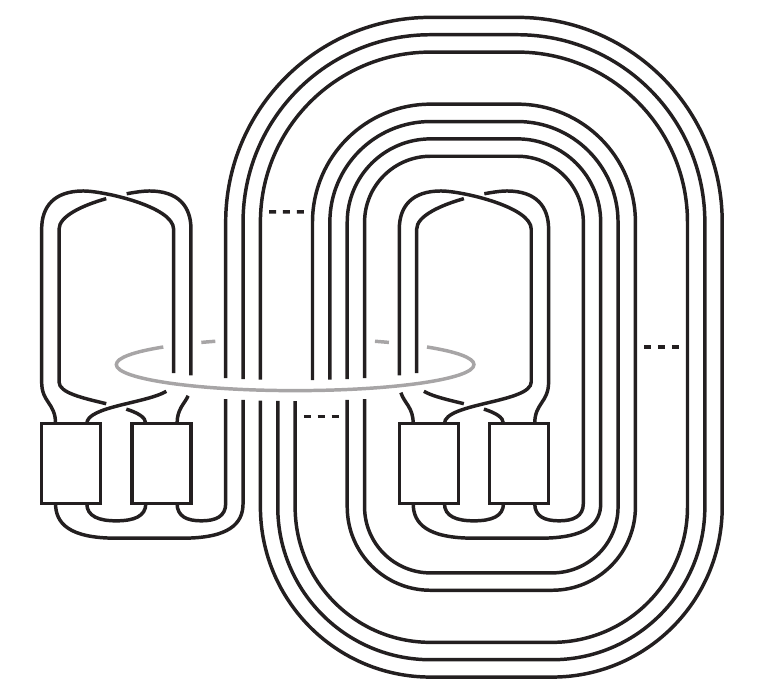}}
\put(8,151){$K_{0,k}$}
\put(30,100){$O$}
\put(62,143){$\overbrace{\hspace{49pt}}^{\colorbox{white}{2k+1}}$}
\put(17,65){$p_1$}
\put(43,65){$p_2$}
\put(120,65){$p_3$}
\put(146,65){$p_4$}
\end{picture}
\caption{
The knots $K_{0,k}$ and $O$.
The sphere $S_{0,k}$ is represented as the horizontal plane containing $O$.}
\label{fig_quotient+}
\end{figure}

The $q$ full twists in the construction of $K_{q,k}$ can be regarded as a $(1/q)$-surgery along a trivial knot $O$ as in Figure~\ref{fig_quotient+}.
That is to say, $(S^3,K_{q,k})$ is obtained from $(S^3,K_{0,k})$ by drilling out a closed regular neighborhood $V$ of $O$, and filling it back so that a meridian of $V$ is identified with a $(1/q)$-curve on $\partial W$, where $W$ denotes the closure of $S^3\setminus V$.
We may suppose that $O$ is contained in the sphere $S_{0,k}$, and that $\partial (S_{0,k}\cap V)$ is identified with $\partial (S_{0,k}\cap W)$ under the surgery.
It follows that the images of $S_{0,k}\cap V$ and $S_{0,k}\cap W$ form the sphere $S_{q,k}$.
Let $V_{q,k}$ and $W_{q,k}$ denote the images of $V$ and $W$, respectively.

\begin{figure}[!ht]
\begin{picture}(210,200)(0,0)
\put(0,0){\includegraphics{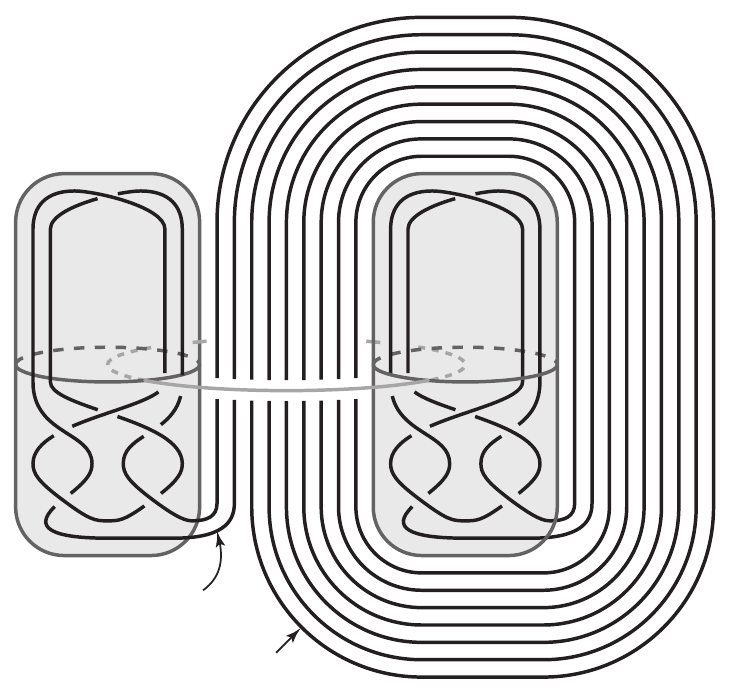}}
\put(26,115){$B_1$}
\put(130,115){$B_2$}
\put(48,25){$a_1$}
\put(70,6){$a_2$}
\end{picture}\\[5pt]
{\Large $\downarrow $}\\[5pt]
\begin{picture}(170,100)(0,0)
\put(0,0){\includegraphics{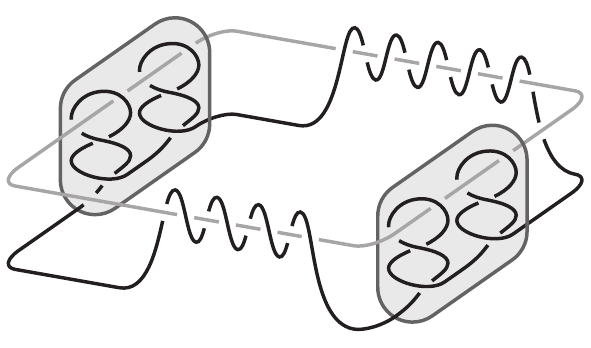}}
\put(24,90){$B_1$}
\put(140,10){$B_2$}
\put(60,15){$a_1$}
\put(126,91){$a_2$}
\end{picture}\\[3pt]
{\Large $\downarrow $}\\[-3pt]
\begin{picture}(340,75)(0,0)
\put(0,0){\includegraphics{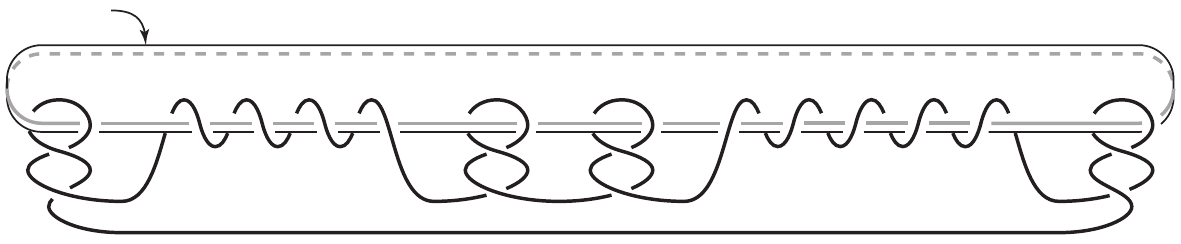}}
\put(12,67){$S_{0,k}$}
\end{picture}\\[5pt]
{\Large $\downarrow $}\\[5pt]
\includegraphics{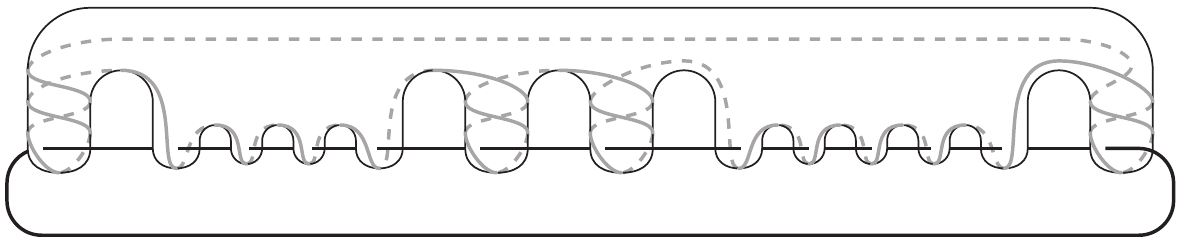}
\caption{
Deformation of the link $K_{0,k}\cup O$ and the sphere $S_{0,k}$, in the case where $p_1=p_2=p_3=p_4=2$ and $k=4$.
In the top two pictures, $S_{0,k}$ is represented as the horizontal plane containing $O$.}
\label{fig_deformation}
\end{figure}

The knot type of $K_{q,k}$ does not depend on the number $k$, which can be seen as follows.
Note that the link $K_{0,k}\cup O$ can be obtained from $K_{0,0}\cup O$ by rotating the ball $B_2$ illustrated in the top of Figure~\ref{fig_deformation} about the axis $O$.
This shows that the link type of $K_{0,k}\cup O$ does not depend on $k$.
The link type and the integer $q$ uniquely determine the knot type of the knot obtained from $K_{0,k}$ by the $(1/q)$-surgery along $O$.
The resultant knot therefore does not depend on $k$.
This, together with the fact mentioned in the last paragraph, gives the desired result.

We deform the link $K_{0,k}\cup O$ together with the sphere $S_{0,k}$ by an ambient isotopy of $S^3$ as illustrated in Figure~\ref{fig_deformation}.
In the first step, the two arcs $a_1,a_2$ of $K_{0,k}\setminus (B_1\cup B_2)$ are pulled forward and backward, respectively, to coil around $O$.
In the second step, the ball above $S_{0,k}$ is shrunk, and $K_{0,k}$ is moved to be suspended from one edge of $O$.
In the third step, $K_{0,k}$ is straightened and $O\subset S_{0,k}$ is distorted for it.
Abusing notation, we continue to denote the results of this deformation by $K_{0,k}$, $S_{0,k}$, $O$, $V$, $W$, and to denote the images under the $(1/q)$-surgery by $K_{q,k}$, $S_{q,k}$, $V_{q,k}$, $W_{q,k}$.

From now on, we consider the $2$-fold covering of $S^3$ branched along $K_{q,k}$, denoted by $\rho _{q,k}:M_{q,k}\rightarrow S^3$.
Let $\Sigma _{q,k}=\rho _{q,k}^{-1}(S_{q,k})$, $\widetilde{V}_{q,k}=\rho _{q,k}^{-1}(V_{q,k})$, $\widetilde{W}_{q,k}=\rho _{q,k}^{-1}(W_{q,k})$.

\begin{figure}[ht]
\begin{picture}(330,100)(0,0)
\put(0,0){\includegraphics{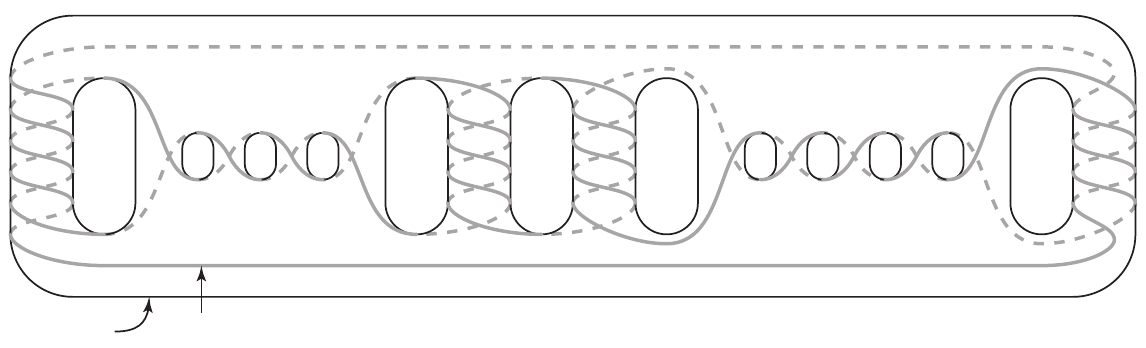}}
\put(13,0){$\Sigma _{0,k}$}
\put(54,0){$P$}
\end{picture}
\caption{The Heegaard surface $\Sigma _{0,k}$ and the pretzel knot $P$, in the case where $p_1=p_2=p_3=p_4=2$ and $k=4$.}
\label{fig_pretzel}
\end{figure}

We can understand the covering $\rho _{0,k}:M_{0,k}\rightarrow S^3$ from the bottom of Figure~\ref{fig_deformation}.
Specifically, we can see the following.
The covering space $M_{0,k}$ is the $3$-sphere, $\Sigma _{0,k}$ is a Heegaard surface of genus $2k+4$, and $\rho _{0,k}^{-1}(O)$ is a knot $P$ contained in $\Sigma _{0,k}$ as illustrated in Figure~\ref{fig_pretzel}.
Note that $P$ is the pretzel knot of type
\begin{equation*}
(2p_1+1,\ \underbrace{1,\ 1,\ \ldots ,\ 1}_k,\ 2p_3+1,\ 2p_4+1,\ \underbrace{-1,\ -1,\ \ldots ,\ -1}_{k+1},\ 2p_2+1).
\end{equation*}
The covering map $\rho _{0,k}$ is the quotient by the $\pi $-rotation of $M_{0,k}$ about the central horizontal axis.
Note that $\widetilde{V}_{0,k}$ is a closed regular neighborhood of $P$ in $M_{0,k}$, and $\widetilde{W}_{0,k}$ is the closure of $M_{0,k}\setminus \widetilde{V}_{0,k}$.

We can regard $\rho _{2q',k}|_{\widetilde{W}_{2q',k}}:\widetilde{W}_{2q',k}\rightarrow W_{2q',k}$ as the same covering as $\rho _{0,k}|_{\widetilde{W}_{0,k}}:\widetilde{W}_{0,k}\rightarrow W_{0,k}$ for each integer $q'$, as follows.
Let $\varphi _{2q',k}:H_1\left( S^3\setminus K_{2q',k}\right) \rightarrow {\mathbb Z}/2{\mathbb Z}$ denote the homomorphism associated with $\rho _{2q',k}$, and $e_{2q',k}:H_1\left( W_{2q',k}\setminus K_{2q',k}\right) \rightarrow H_1\left( S^3\setminus K_{2q',k}\right) $ denote the natural epimorphism.
Since $W_{2q',k}\setminus K_{2q',k}$ coincides with $W_{0,k}\setminus K_{0,k}$, we may regard $e_{2q',k}$ as defined on $H_1\left( W_{0,k}\setminus K_{0,k}\right) $.
Note that $H_1\left( W_{0,k}\setminus K_{0,k}\right) $ is generated by the elements coming from meridians of $K_{0,k}$ and $O$, denoted by $m$ and $\mu $, respectively.
Let $\lambda $ denote the element coming from a longitude of $O$.
We have $(\varphi _{2q',k}\circ e_{2q',k})(\mu )=(\varphi _{2q',k}\circ e_{2q',k})(\mu )+2(\varphi _{2q',k}\circ e_{2q',k})(q'\lambda )=(\varphi _{2q',k}\circ e_{2q',k})(\mu +2q'\lambda )=0$ since $e_{2q',k}(\mu +2q'\lambda )$ comes from a meridian of $V_{2q',k}$.
Also, we have $(\varphi _{2q',k}\circ e_{2q',k})(m)=1$ since $e_{2q',k}(m)$ comes from a meridian of the branch set $K_{2q',k}$ of $\rho _{2q',k}$.
Thus, we obtain $\varphi _{2q',k}\circ e_{2q',k}=\varphi _{0,k}\circ e_{0,k}$.
This shows that the associated coverings $\rho _{2q',k}|_{\widetilde{W}_{2q',k}}$ and $\rho _{0,k}|_{\widetilde{W}_{0,k}}$ coincide under a homeomorphism between $\widetilde{W}_{2q',k}$ and $\widetilde{W}_{0,k}$.

We can regard $\rho _{2q',k}|_{\widetilde{V}_{2q',k}}:\widetilde{V}_{2q',k}\rightarrow V_{2q',k}$ as the same covering as $\rho _{0,k}|_{\widetilde{V}_{0,k}}:\widetilde{V}_{0,k}\rightarrow V_{0,k}$, as follows.
Since $V_{2q',k}$ is disjoint from $K_{2q',k}$, the restriction $\rho _{2q',k}|_{\widetilde{V}_{2q',k}}$ is an unbranched $2$-fold covering of the solid torus $V_{2q',k}$.
There are two possibilities in one of which $\widetilde{V}_{2q',k}$ is a solid torus, and in the other of which $\widetilde{V}_{2q',k}$ is the union of two solid tori.
Note that $\partial \widetilde{W}_{2q',k}$ is a single torus as well as $\partial \widetilde{W}_{0,k}$ by the above results.
Since $\partial \widetilde{V}_{2q',k}$ is identified with the single torus, $\widetilde{V}_{2q',k}$ is a single solid torus as well as $\widetilde{V}_{0,k}$.

The total covering space $M_{2q',k}$ is obtained by gluing $\widetilde{W}_{2q',k}$ and $\widetilde{V}_{2q',k}$ appropriately.
Recall that a meridian of $V_{2q',k}$ is identified with a $(1/2q')$-curve on $\partial W_{2q',k}$, and $\partial (S_{2q',k}\cap V_{2q',k})$ is identified with $\partial (S_{2q',k}\cap W_{2q',k})$.
Note that the meridian of $V_{2q',k}$ lifts to a pair of meridians of $\widetilde{V}_{2q',k}$.
One can see that a $(1/2q')$-curve on $\partial W_{0,k}$ lifts to a pair of $(1/q')$-curves on $\partial \widetilde{W}_{0,k}$, and it follows that the $(1/2q')$-curve on $\partial W_{2q',k}$ lifts to a pair of $(1/q')$-curves on $\partial \widetilde{W}_{2q',k}$.
These facts show that a meridian of $\widetilde{V}_{2q',k}$ is identified with a $(1/q')$-curve on $\partial \widetilde{W}_{2q',k}$, and $\partial (\Sigma _{2q',k}\cap \widetilde{V}_{2q',k})$ is identified with $\partial (\Sigma _{2q',k}\cap \widetilde{W}_{2q',k})$.

By the above, the covering space $M_{2q',k}$ and the lift $\Sigma _{2q',k}$ of $S_{2q',k}$ turn out to be the $3$-manifold and the Heegaard surface constructed by Casson--Gordon.
That is to say, $M_{2q',k}$ is obtained from the $3$-sphere $M_{0,k}$ by a $(1/q')$-surgery along the pretzel knot $P$, and $\Sigma _{2q',k}$ is the image of the Heegaard surface $\Sigma _{0,k}$.
Casson--Gordon showed that if $|2p_i+1|\geq 5$ for each $i\in \{ 1,2,3,4\} $ and $|q'|\geq 6$, then the Heegaard surface $\Sigma _{2q',k}$ is strongly irreducible.
Though their original proof has been unpublished, one can verify this fact by \cite[Theorem A]{Moriah-Schultens} and the incompressibility of $\Sigma _{0,k}\setminus P$ shown in \cite{Paris}.

The strong irreducibility of $\Sigma _{q,k}$ guarantees the strong irreducibility of $S_{q,k}$, which can be seen as follows.
Assume that $S_{q,k}$ is not strongly irreducible.
There exist disjoint essential disks $D^+,D^-$ of the trivial tangles divided by $S_{q,k}$.
A component of $\rho _{q,k}^{-1}(D^+)$ and a component of $\rho _{q,k}^{-1}(D^-)$ are disjoint essential disks of the handlebodies divided by $\Sigma _{q,k}$.
This contradicts the strong irreducibility of $\Sigma _{q,k}$.

Thus, we conclude that if $|2p_i+1|\geq 5$ for each $i\in \{ 1,2,3,4\} $ and $q$ is an even number with $|q|\geq 12$, then the bridge sphere $S_{q,k}$ is strongly irreducible, and hence destabilized, for $k\geq 0$.
This completes the proof of Theorem~\ref{quotient}.

We remark that $\{K_{2q',k}\}_{q'\in\mathbb{Z}}$ contains infinitely many equivalence classes of knots, which can be seen as follows.
By \cite{Wu}, the knot $P$ is hyperbolic and the manifold $M_{2q',k}$ is hyperbolic if $|2p_i+1|\geq 2$ for each $i\in \{ 1,2,3,4\} $ and $q'\ne 0$.
By \cite{Neumann-Zagier}, the set $\{ \text{the length of the shortest geodesic in }M_{2q',k}\} _{q'\in\mathbb{Z}}$ is an infinite set.

We remark that another proof of Theorem~\ref{quotient} by constructing our bridge spheres from Casson--Gordon's Heegaard surfaces is no easier than the above proof.
It is well known that one can construct a bridge sphere of a link from a Heegaard surface of a $3$-manifold if the manifold admits an involution whose restriction to each of the handlebodies is hyperelliptic.
It is true that the $\pi $-rotation of $W_{2q',k}=W_{0,k}\subset M_{0,k}$ about the central horizontal axis in Figure~\ref{fig_pretzel} extends to an involution of $M_{2q',k}$.
It is also true, but cannot be proved easily that the restriction of the involution to each of the handlebodies divided by $\Sigma _{2q',k}$ is hyperelliptic.

\section{Bridge diagrams and $2$-connected condition}\label{diagrams}

In this section, we define the notions of bridge diagram and $2$-connected condition, and prove the following theorem.

\begin{Thm}\label{criterion}
An $n$-bridge sphere $S$ with $n\geq 3$ of a link is strongly irreducible if there exist a bridge diagram $(\sigma ^+,\sigma ^-)$ of the bridge sphere $S$ and a loop $l$ on $S$ such that the pair $\left( (\sigma ^+,\sigma ^-),l\right) $ satisfies the $2$-connected condition.
\end{Thm}

\noindent
This is a refinement of previous work by the fourth author \cite{Takao}.
We mention that Kwon \cite{Kwon} recently gave another criterion for recognizing strongly irreducible bridge spheres.

We use the following notation.
Let $S$ be an $n$-bridge sphere with $n\geq 3$ of a link $L$ in $S^3$.
Let $B^+,B^-\subset S^3$ denote the $3$-balls divided by $S$, and let $\tau ^+=L\cap B^+$, $\tau ^-=L\cap B^-$.
For each $\varepsilon \in \{ +,-\} $, let $\{ E_1^\varepsilon ,E_2^\varepsilon ,\ldots ,E_n^\varepsilon \} $ be a complete collection of bridge disks of $(B^\varepsilon ,\tau ^\varepsilon )$, and let $\tau _i^\varepsilon =E_i^\varepsilon \cap L$, $\sigma _i^\varepsilon =E_i^\varepsilon \cap S$, and $\sigma ^\varepsilon =\sigma _1^\varepsilon \cup \sigma _2^\varepsilon \cup \cdots \cup \sigma _n^\varepsilon $.

A {\it bridge diagram} of the bridge sphere $S$ is the diagram composed of the arcs of $\sigma ^+$ and $\sigma ^-$ on $S$, and is denoted by $(\sigma ^+,\sigma ^-)$.
One can think of it as a link diagram of $L$ obtained by projecting $\tau _1^\varepsilon ,\tau _2^\varepsilon ,\ldots ,\tau _n^\varepsilon $ into $S$ disjointly for each $\varepsilon \in \{ +,-\} $.

Let $l$ be a simple loop on $S$ containing $\sigma ^-$ such that the arcs $\sigma _1^-,\sigma _2^-,\ldots ,\sigma _n^-$ are located in $l$ in this order.
We may suppose that $\sigma ^+$ intersects $l$ transversely and $|\sigma ^+\cap l|$ is minimal by deforming $E_1^+,E_2^+,\ldots ,E_n^+$ by an ambient isotopy of $B^+$ rel $L\cap B^+$.

We define a graph ${\mathcal G}_{i,j,\varepsilon }$ for distinct $i,j\in \{ 1,2,\ldots ,n\} $ and $\varepsilon \in \{ +,-\} $ as follows.
Let $S_+,S_-\subset S$ denote the hemispheres divided by $l$, and $\delta _i$ denote the component of $l\setminus \sigma ^-$ between $\sigma _i^-$ and $\sigma _{i+1}^-$ for each $i\in \{ 1,2,\ldots n\} $, where the index $i+1$ is considered modulo $n$.
Let ${\mathcal A}_{i,j,\varepsilon }$ be the collection of components of $\sigma ^+\cap S_\varepsilon $ separating $\delta _i$ from $\delta _j$ in $S_\varepsilon $ for distinct $i,j\in \{ 1,2,\ldots ,n\} $ and $\varepsilon \in \{ +,-\} $.
Note that ${\mathcal A}_{i,j,\varepsilon }$ consists of parallel arcs in $S_\varepsilon $.
Let ${\mathcal G}_{i,j,\varepsilon }$ be the graph such that the vertex set is $\{ 1,2,\ldots ,n\} $, and distinct vertices $v,w$ span an edge if ${\mathcal A}_{i,j,\varepsilon }$ has a subarc of $\sigma _v^+$ and a subarc of $\sigma _w^+$ that are adjacent in $S_\varepsilon $.
By ``adjacent'', we mean that the subarcs of $\sigma _v^+$ and $\sigma _w^+$ are contained in the closure of one component of $S_\varepsilon \setminus \sigma ^+$.

We say that the pair $\left( (\sigma ^+,\sigma ^-),l\right) $ satisfies the {\it $2$-connected condition} if the graph ${\mathcal G}_{i,j,\varepsilon }$ is $2$-connected for every combination of distinct $i,j\in \{ 1,2,\ldots ,n\} $ and $\varepsilon \in \{ +,-\} $.
Here, a graph is said to be $2$-connected if it is connected even after deleting any vertex.

\begin{figure}[ht]
\begin{picture}(280,100)(0,0)
\put(0,0){\includegraphics{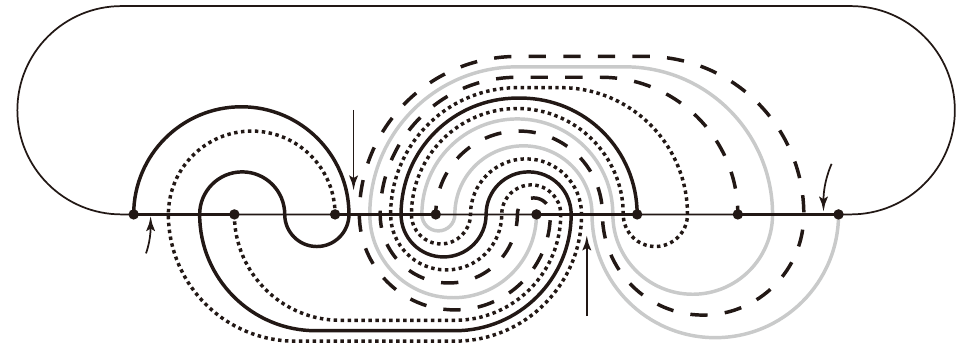}}
\put(15,50){$l$}
\put(20,75){$S_+$}
\put(5,20){$S_-$}
\put(62,75){$\sigma _1^+$}
\put(45,3){$\sigma _2^+$}
\put(220,75){$\sigma _3^+$}
\put(231,3){$\sigma _4^+$}
\put(35,18){$\sigma _1^-$}
\put(98,75){$\sigma _2^-$}
\put(164,0){$\sigma _3^-$}
\put(238,58){$\sigma _4^-$}
\end{picture}\\[15pt]
\begin{picture}(280,75)(0,0)
\put(0,5){\includegraphics{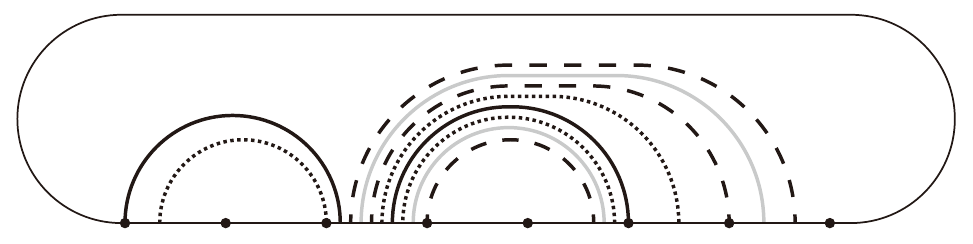}}
\put(20,50){$S_+$}
\put(44,0){$\sigma _1^-$}
\put(75,0){$\delta _1$}
\put(102,0){$\sigma _2^-$}
\put(133,0){$\delta _2$}
\put(162,0){$\sigma _3^-$}
\put(193,0){$\delta _3$}
\put(220,0){$\sigma _4^-$}
\put(265,8){$\delta _4$}
\end{picture}\\[15pt]
\begin{picture}(65,65)(0,0)
\put(0,0){\includegraphics{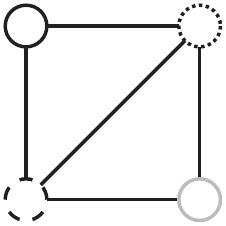}}
\put(5,54){$1$}
\put(55,54){$2$}
\put(5,4){$3$}
\put(55,4){$4$}
\end{picture}
\caption{A bridge diagram, the arcs of ${\mathcal A}_{1,2,+}$ and the graph ${\mathcal G}_{1,2,+}$.}
\label{fig_criterion}
\end{figure}

For example, we examine the bridge diagram $(\sigma ^+,\sigma ^-)$ of a $4$-bridge sphere of a knot as in the top of Figure~\ref{fig_criterion}.
We choose $l$, $S_+$, $S_-$ as in the figure.
The middle of the figure shows the arcs of ${\mathcal A}_{1,2,+}$, and the bottom shows the graph ${\mathcal G}_{1,2,+}$.
Clearly, ${\mathcal G}_{1,2,+}$ is $2$-connected.
One can check that ${\mathcal G}_{2,3,+}$ is, however, not $2$-connected, and so the pair $\left( (\sigma ^+,\sigma ^-),l\right) $ does not satisfy the $2$-connected condition.
See Section~\ref{proof_biplat} for examples satisfying the $2$-connected condition.

In the rest of this section, we give a proof of Theorem~\ref{criterion}.

By way of contradiction, we assume that an $n$-bridge sphere $S$ is not strongly irreducible and the pair $\left( (\sigma ^+,\sigma ^-),l\right) $ satisfies the $2$-connected condition with the above notation.

There exist disjoint essential disks $D^+$, $D^-$ of the tangles $(B^+,\tau ^+)$, $(B^-,\tau ^-)$, respectively, since $S$ is not strongly irreducible.
By \cite[Lemma~3]{Takao}, we may simultaneously suppose that
\begin{gather*}
|\partial D^+\cap\sigma ^+|={\rm min}\left\{ |\partial \widetilde{D}^+\cap\sigma ^+| \ \right| \left. \text{$\widetilde{D}^+$ is isotopic to $D^+$ in $B^+$ rel $L\cap B^+$}\right\} ,\\
|\partial D^-\cap l|={\rm min}\left\{ |\partial \widetilde{D}^-\cap l| \ \right| \left. \text{$\widetilde{D}^-$ is isotopic to $D^-$ in $B^-$ rel $L\cap B^-$}\right\} .
\end{gather*}

We consider another graph ${\mathcal G}(D^-)$ defined as follows.
The vertex set of ${\mathcal G}(D^-)$ is $\{ 1,2,\ldots ,n\} $.
Distinct vertices $v,w$ of ${\mathcal G}(D^-)$ span an edge if $\partial D^-$ contains a subarc connecting $\sigma _v^+$ and $\sigma _w^+$ directly.
By ``directly'', we mean that the interior of the subarc is disjoint from $\sigma ^+$.

On the one hand, we show that the graph ${\mathcal G}(D^-)$ is $2$-connected as follows.
By \cite[Lemma~1]{Takao}, there exist distinct $i,j\in \{ 1,2,\ldots ,n\} $ and $\varepsilon \in \{ +,-\} $ such that $\partial D^-$ contains a subarc connecting $\delta _i$ and $\delta _j$ in $S_\varepsilon $.
Note that the subarc of $\partial D^-$ intersects all the arcs of ${\mathcal A}_{i,j,\varepsilon }$.
In particular, if two arcs of ${\mathcal A}_{i,j,\varepsilon }$ are adjacent in $S_\varepsilon $, then $\partial D^-$ contains a subarc connecting them directly.
This shows that ${\mathcal G}_{i,j,\varepsilon }$ is a subgraph of ${\mathcal G}(D^-)$.
Recall that ${\mathcal G}_{i,j,\varepsilon }$ and ${\mathcal G}(D^-)$ have the same vertex set $\{ 1,2,\ldots ,n\} $.
The graph ${\mathcal G}(D^-)$ therefore inherits the $2$-connectedness from ${\mathcal G}_{i,j,\varepsilon }$.

On the other hand, we show that the graph ${\mathcal G}(D^-)$ is not $2$-connected as follows.
After an isotopy of ${\rm int}D^+$ if necessary, we may suppose that $D^+$ intersects $E_1^+\cup E_2^+\cup \cdots \cup E_n^+$ transversely and $D^+\cap (E_1^+\cup E_2^+\cup \cdots \cup E_n^+)$ has no loop component.
First, consider the case where $D^+\cap (E_1^+\cup E_2^+\cup \cdots \cup E_n^+)$ is empty.
Note that each of $E_1^+,E_2^+,\ldots ,E_n^+$ is contained in one of the two components of $B^+\setminus D^+$.
This gives a decomposition of $\{ 1,2,\ldots ,n\} $ into two subsets $V$ and $W$.
Since $D^+$ is an essential disk of $(B^+,\tau ^+)$, both $V$ and $W$ are not empty.
By the disjointness of $D^+$ and $D^-$, for any $v\in V$ and $w\in W$, the curve $\partial D^-$ cannot connect $\sigma _v^+$ and $\sigma _w^+$.
This implies that the graph ${\mathcal G}(D^-)$ is not connected, and hence not 2-connected.
Secondly, consider the case where $D^+\cap (E_1^+\cup E_2^+\cup \cdots \cup E_n^+)$ is not empty.
Note that it consists of properly embedded arcs in $D^+$.
Let $D_0^+$ be an outermost subdisk of $D^+$ cut off by an arc of $D^+\cap (E_1^+\cup E_2^+\cup \cdots \cup E_n^+)$, which is an arc of $D^+\cap E_u^+$ for some $u\in \{ 1,2,\ldots ,n\} $.
Note that each of $E_1^+,\ldots ,E_{u-1}^+,E_{u+1}^+,\ldots ,E_n^+$ is contained in one of the two components of $B^+\setminus (D_0^+\cup E_u^+)$.
This gives a decomposition of $\{ 1,\ldots ,u-1,u+1,\ldots ,n\} $ into two subsets $V$ and $W$.
Since $n\geq 3$ and $|\partial D^+\cap \sigma ^+|$ is minimal, both $V$ and $W$ are not empty.
By the disjointness of $D^+$ and $D^-$, for any $v\in V$ and $w\in W$, any subarc of $\partial D^-$ cannot connect $\sigma _v^+$ and $\sigma _w^+$ directly.
This implies that the graph ${\mathcal G}(D^-)$ is not connected after deleting the vertex $u$, and hence ${\mathcal G}(D^-)$ is not 2-connected.

Thus, we have a contradiction to the $2$-connectedness of the graph ${\mathcal G}(D^-)$.
This completes the proof of Theorem~\ref{criterion}.

\section{Proof of Theorem~2}\label{proof_biplat}

In this section, we give a proof of Theorem~\ref{biplat}.
To do this, we use the method introduced in Section~\ref{diagrams}.
This is a sequel of previous work by the third and fourth authors \cite{Ozawa-Takao}.

First of all, it is easy to see that the spheres $S_n$ and $S_n^\perp $ shown in Figure~\ref{fig_biplat} are an $n$-bridge sphere and a $3$-bridge sphere, respectively, of the knot $K_n$.
Indeed, the figure shows a $2n$-plat presentation of $K_n$, and shows a $6$-plat presentation of $K_n$ after the $(\pi /2)$-rotation and a slight deformation (see \cite{Birman1} for the concept of plat presentation).

\begin{figure}[ht]
\begin{picture}(290,175)(0,0)
\put(0,10){\includegraphics{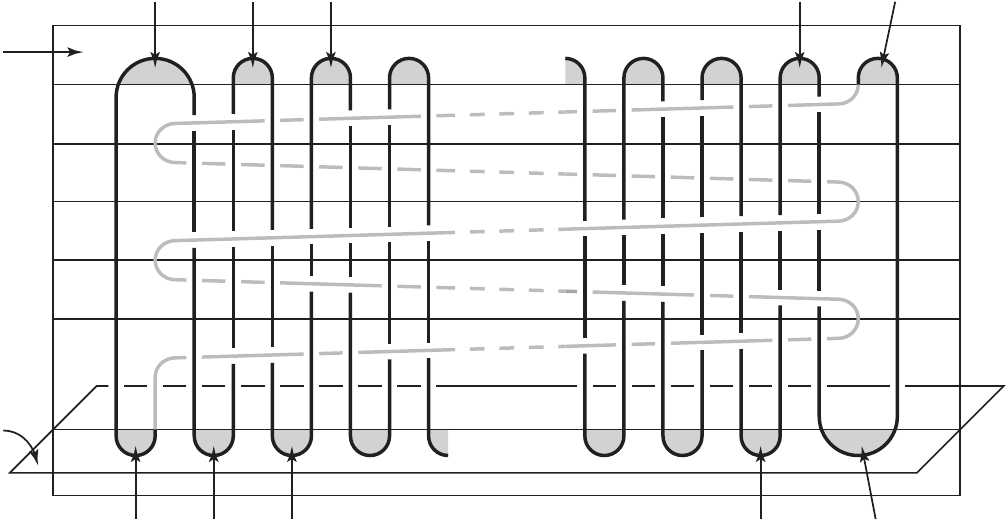}}
\put(-11,143){$A$}
\put(-13,34){$S_n$}
\put(280,150){$6$}
\put(280,133){$5$}
\put(280,115){$4$}
\put(280,99){$3$}
\put(280,82){$2$}
\put(280,65){$1$}
\put(281,31){$0$}
\put(280,13){$-1$}
\put(40,165){$\bar{E}_1^+$}
\put(68,165){$\bar{E}_2^+$}
\put(90,165){$\bar{E}_3^+$}
\put(109,165){$\cdots $}
\put(226,165){$\bar{E}_{n-1}^+$}
\put(255,165){$\bar{E}_n^+$}
\put(33,0){$E_1^-$}
\put(56,0){$E_2^-$}
\put(79,0){$E_3^-$}
\put(98,0){$\cdots $}
\put(214,0){$E_{n-1}^-$}
\put(248,0){$E_n^-$}
\end{picture}\\[3pt]
in the case where $n$ is even\\[12pt]
\begin{picture}(290,175)(0,0)
\put(0,10){\includegraphics{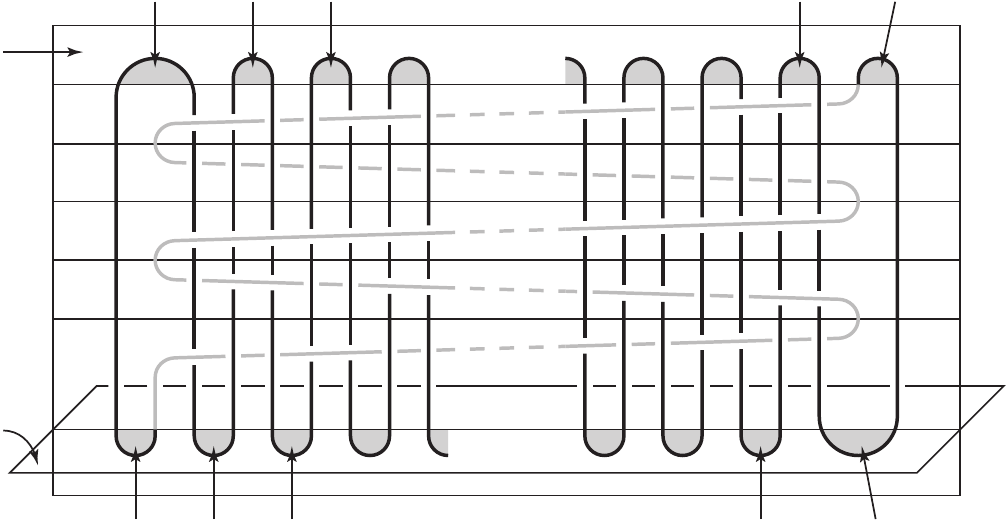}}
\put(-11,143){$A$}
\put(-13,34){$S_n$}
\put(280,150){$6$}
\put(280,133){$5$}
\put(280,115){$4$}
\put(280,99){$3$}
\put(280,82){$2$}
\put(280,65){$1$}
\put(281,31){$0$}
\put(280,13){$-1$}
\put(40,165){$\bar{E}_1^+$}
\put(68,165){$\bar{E}_2^+$}
\put(90,165){$\bar{E}_3^+$}
\put(109,165){$\cdots $}
\put(226,165){$\bar{E}_{n-1}^+$}
\put(255,165){$\bar{E}_n^+$}
\put(33,0){$E_1^-$}
\put(56,0){$E_2^-$}
\put(79,0){$E_3^-$}
\put(98,0){$\cdots $}
\put(214,0){$E_{n-1}^-$}
\put(248,0){$E_n^-$}
\end{picture}\\[3pt]
in the case where $n$ is odd
\caption{The knot $K_n$ and a vertical annulus $A$ in the product space $S_n\times [-1,6]$, and disks $\bar{E}_1^+,\bar{E}_2^+,\ldots ,\bar{E}_n^+$, $E_1^-,E_2^-,\ldots ,E_n^-$.}
\label{fig_biplat+}
\end{figure}

We regard the knot $K_n$ as lying in the product space $S_n\times [-1,6]$ as in Figure~\ref{fig_biplat+}.
The sphere $S_n$ is regarded as the section $S_n\times \{ 0\} $, and we let $S_n(s)=S_n\times \{ s\} $ for each $s\in [0,5]$.
We suppose that $K_n$ is contained in the vertical annulus $A$ shown in the figure except for the gray arc.
Let $l(s)=A\cap S_n(s)$ for each $s\in [0,5]$ and let $l=l(0)$.

We choose complete collections of bridge disks of the trivial tangles divided by $S_n$ as follows.
Let $\bar{E}_1^+,\bar{E}_2^+,\ldots ,\bar{E}_n^+$ and $E_1^-,E_2^-,\ldots ,E_n^-$ be the disks in the vertical annulus $A$ as in Figure~\ref{fig_biplat+}.
Note that $\{ E_1^-,E_2^-,\ldots ,E_n^-\} $ is a complete collection of bridge disks of the trivial tangle below $S_n$.
Let $\sigma _i^-=E_i^-\cap S_n$ for each $i\in \{ 1,2,\ldots ,n\} $ and $\sigma ^-=\sigma _1^-\cup \sigma _2^-\cup \cdots \cup \sigma _n^-$.
There exists a homeomorphism $\Phi :S_n\times [0,5]\rightarrow S_n\times [0,5]$ such that $\Phi \left(S_n(s)\right) =S_n(s)$ for $s\in [0,5]$, $\Phi |_{S_n(5)}=id_{S_n(5)}$ and $\Phi \left(K_n\cap \left( S_n\times [0,5]\right) \right) =\left( K_n\cap S_n(5)\right) \times [0,5]$.
Let $\sigma _i^+(s)$ denote the arc $\Phi ^{-1}\left( \left( \bar{E}_i^+\cap S_n(5)\right) \times \{ s\} \right) $ for each $i\in \{ 1,2,\ldots ,n\} $ and $s\in [0,5]$, and let $\sigma ^+(s)=\sigma _1^+(s)\cup \sigma _2^+(s)\cup \cdots \cup \sigma _n^+(s)$.
We may suppose that $\sigma ^+(r)$ intersects $l(r)$ transversely and $|\sigma ^+(r)\cap l(r)|$ is minimal for each $r\in \{ 0,1,2,3,4\} $ after a modification of $\Phi $.
Let $E_i^+$ denote the disk $\bar{E}_i^+\cup \bigcup _{s\in [0,5]}\sigma _i^+(s)$ for each $i\in \{ 1,2,\ldots ,n\} $.
Note that $\{ E_1^+,E_2^+,\ldots ,E_n^+\} $ is a complete collection of bridge disks of the trivial tangle above $S_n$.
Let $\sigma _i^+=\sigma _i^+(0)$ for each $i\in \{ 1,2,\ldots ,n\} $ and $\sigma ^+=\sigma ^+(0)$.
Note that $\sigma _i^+=E_i^+\cap S_n$.

\begin{figure}[ht]
\begin{picture}(340,40)(0,0)
\put(17,0){\includegraphics{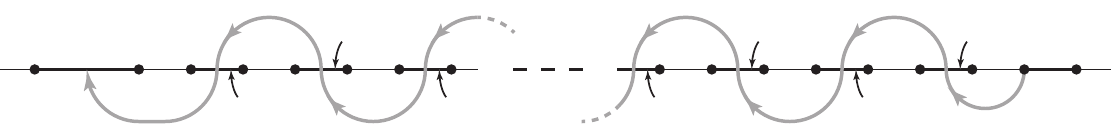}}
\put(0,17){$l(5)$}
\put(47,23){$1$}
\put(85,3){$2$}
\put(116,30){$3$}
\put(146,3){$4$}
\put(201,3){$n-4$}
\put(232,30){$n-3$}
\put(260,3){$n-2$}
\put(291,30){$n-1$}
\put(318,11){$n$}
\end{picture}\\
in the case where $n$ is even\\[15pt]
\begin{picture}(340,40)(0,0)
\put(17,0){\includegraphics{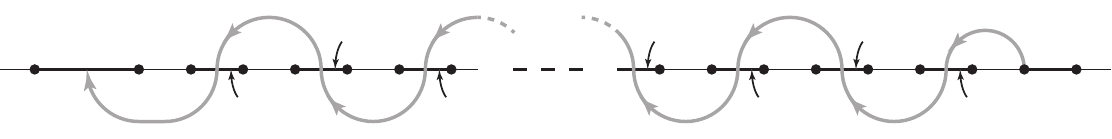}}
\put(0,17){$l(5)$}
\put(47,23){$1$}
\put(85,3){$2$}
\put(116,30){$3$}
\put(146,3){$4$}
\put(201,30){$n-4$}
\put(232,3){$n-3$}
\put(260,30){$n-2$}
\put(291,3){$n-1$}
\put(318,23){$n$}
\end{picture}\\
in the case where $n$ is odd
\caption{The arcs of $\sigma ^+(5)$ in the sphere $S_n(5)$, where ``$i$'' stands for $\sigma _i^+(5)$ for each $i\in \{ 1,2,\ldots ,n\} $.}
\label{fig_diagram_0}
\end{figure}

\begin{figure}[p]
\begin{picture}(290,100)(0,0)
\put(17,0){\includegraphics{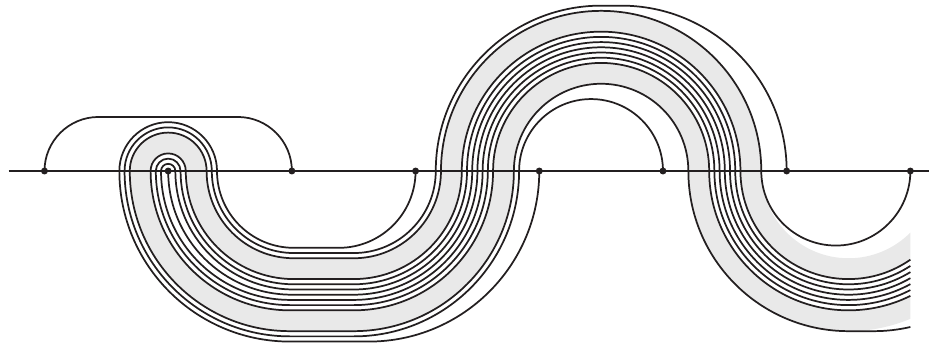}}
\put(0,48){$l(4)$}
\put(98,68){$1$}
\put(174,27){$2$}
\put(247,68){$3$}
\put(264,37){$4$}
\end{picture}\\[5pt]
the leftmost part\\[12pt]
\begin{picture}(290,100)(0,0)
\put(0,0){\includegraphics{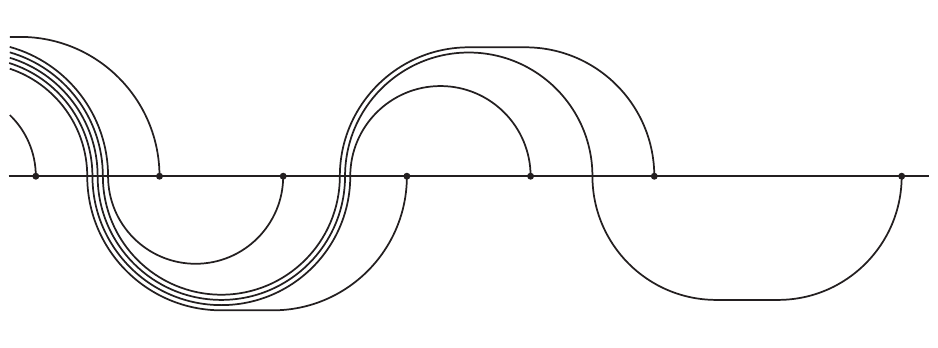}}
\put(272,46){$l(4)$}
\put(47,68){$n-3$}
\put(119,27){$n-2$}
\put(191,68){$n-1$}
\put(245,27){$n$}
\end{picture}\\[5pt]
the rightmost part in the case where $n$ is even\\[12pt]
\begin{picture}(290,100)(0,0)
\put(0,0){\includegraphics{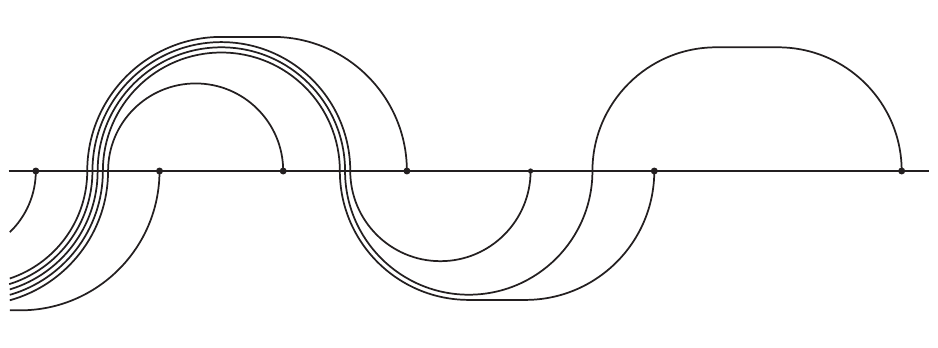}}
\put(272,47){$l(4)$}
\put(47,27){$n-3$}
\put(119,68){$n-2$}
\put(191,27){$n-1$}
\put(245,68){$n$}
\end{picture}\\[5pt]
the rightmost part in the case where $n$ is odd
\caption{
The arcs of $\sigma ^+(4)$ in the sphere $S_n(4)$.
The gray band represents a bunch of parallel arcs.}
\label{fig_diagram_1}
\end{figure}

\begin{figure}[p]
\begin{picture}(290,130)(0,0)
\put(17,0){\includegraphics{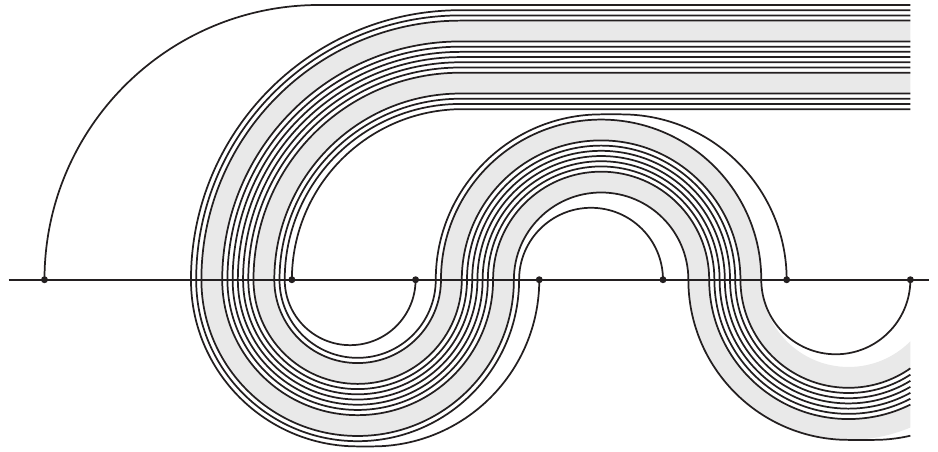}}
\put(0,47){$l(3)$}
\end{picture}\\[5pt]
the leftmost part\\[12pt]
\begin{picture}(290,130)(0,0)
\put(0,0){\includegraphics{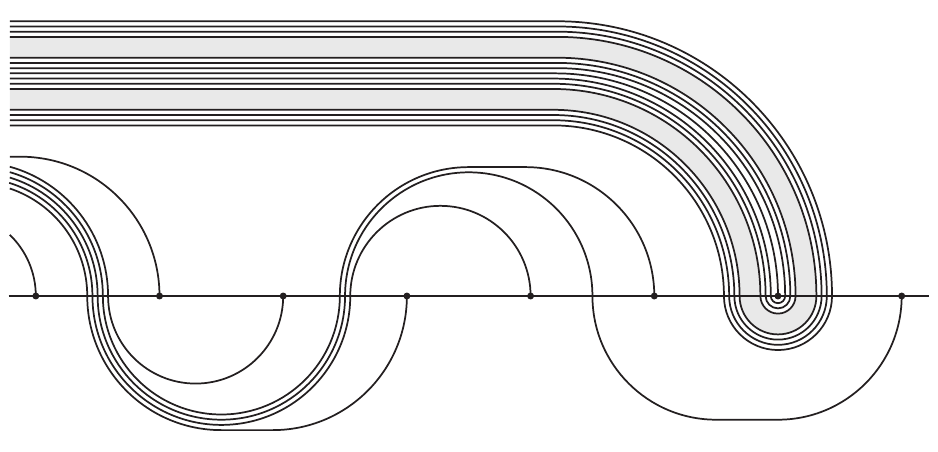}}
\put(272,42){$l(3)$}
\end{picture}\\[5pt]
the rightmost part in the case where $n$ is even\\[12pt]
\begin{picture}(290,130)(0,0)
\put(0,0){\includegraphics{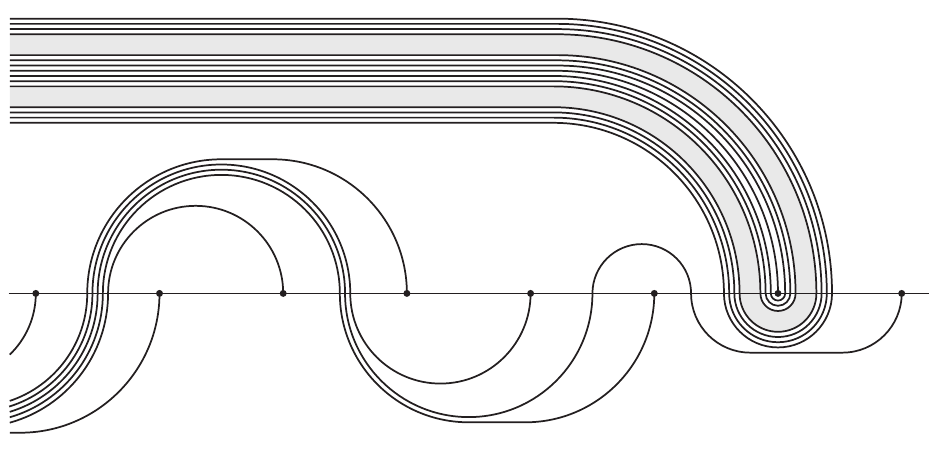}}
\put(272,42){$l(3)$}
\end{picture}\\[5pt]
the rightmost part in the case where $n$ is odd
\caption{The arcs of $\sigma ^+(3)$ in the sphere $S_n(3)$.}
\label{fig_diagram_2}
\end{figure}

\begin{figure}[p]
\begin{picture}(290,140)(0,0)
\put(17,0){\includegraphics{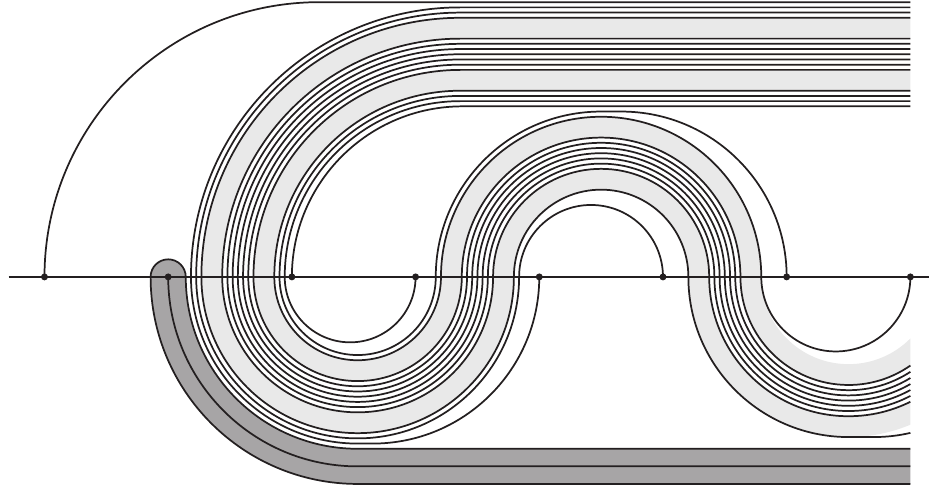}}
\put(0,57){$l(2)$}
\end{picture}\\[5pt]
the leftmost part\\[12pt]
\begin{picture}(290,140)(0,0)
\put(0,0){\includegraphics{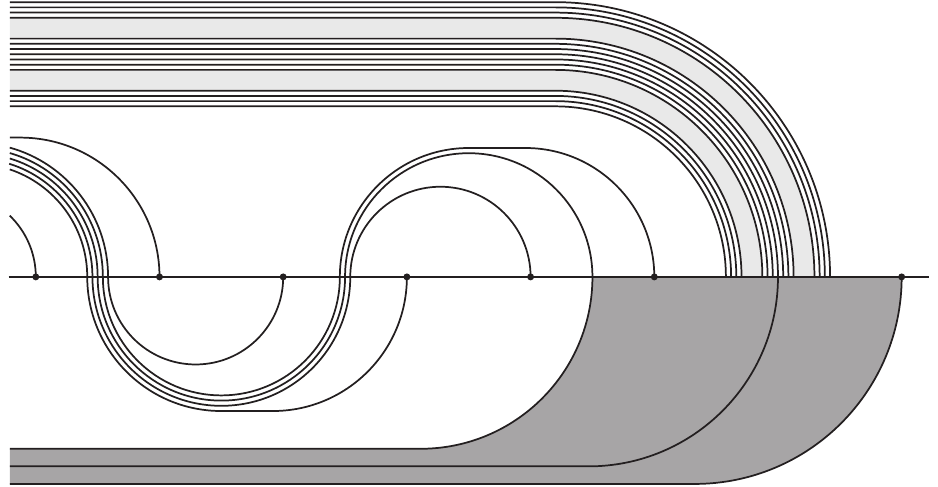}}
\put(272,57){$l(2)$}
\put(230,40){$R(2)$}
\end{picture}\\[5pt]
the rightmost part in the case where $n$ is even\\[12pt]
\begin{picture}(290,140)(0,0)
\put(0,0){\includegraphics{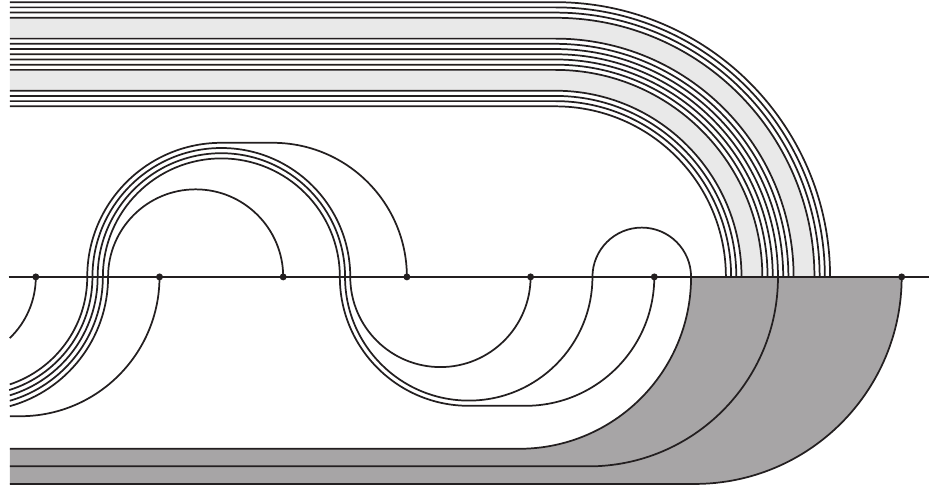}}
\put(272,57){$l(2)$}
\put(230,40){$R(2)$}
\end{picture}\\[5pt]
the rightmost part in the case where $n$ is odd
\caption{
The arcs of $\sigma ^+(2)$ in the sphere $S_n(2)$.
A region $R(2)$ is painted over with darker gray.}
\label{fig_diagram_3}
\end{figure}

\begin{figure}[p]
\begin{picture}(230,110)(0,0)
\put(0,0){\includegraphics{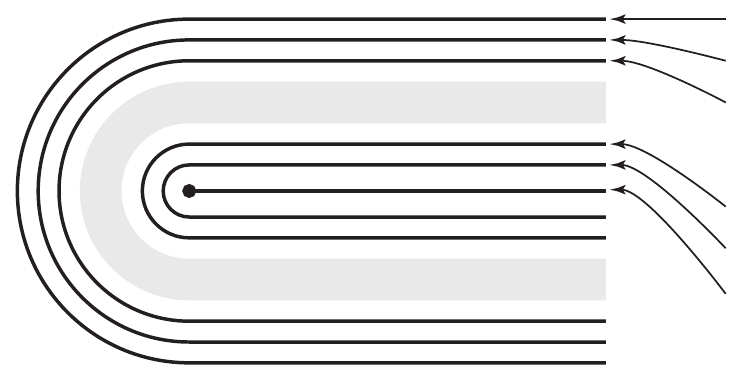}}
\put(212,103){$n$}
\put(212,90){$1$}
\put(212,77){$2$}
\put(214,60){$\vdots $}
\put(212,48){$n-2$}
\put(212,35){$n-1$}
\put(212,22){$n$}
\end{picture}
\caption{Subarcs of $\sigma ^+(r)$ contained in the region $R(r)$ for each $r\in \{ 2,1,0\} $, where ``$i$'' stands for a subarc of $\sigma _i^+(r)$ for each $i\in \{ 1,2,\ldots ,n\} $.}
\label{fig_tongue}
\end{figure}

\begin{figure}[p]
\begin{picture}(290,105)(0,0)
\put(17,0){\includegraphics{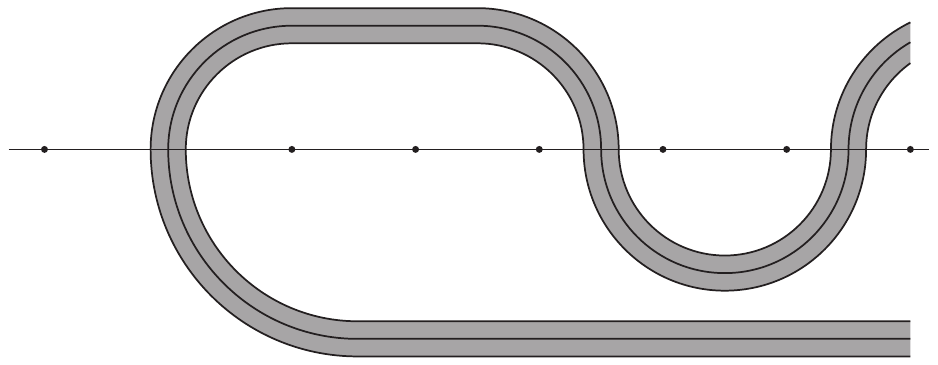}}
\put(0,59){$l(1)$}
\end{picture}\\[5pt]
the leftmost part in the case where $n$ is even\\[12pt]
\begin{picture}(290,105)(0,0)
\put(17,0){\includegraphics{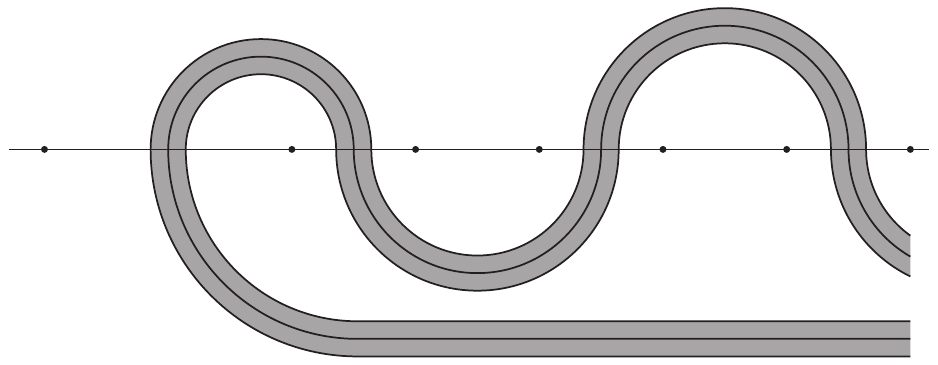}}
\put(0,59){$l(1)$}
\end{picture}\\[5pt]
the leftmost part in the case where $n$ is odd\\[12pt]
\begin{picture}(290,105)(0,0)
\put(0,0){\includegraphics{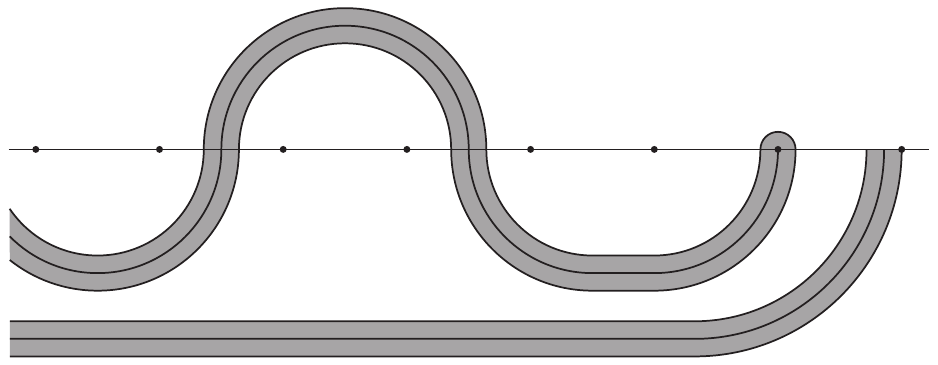}}
\put(272,59){$l(1)$}
\end{picture}\\[5pt]
the rightmost part
\caption{The region $R(1)$ in the sphere $S_n(1)$.}
\label{fig_diagram_4}
\end{figure}

\begin{figure}[p]
\begin{picture}(280,100)(0,0)
\put(7,0){\includegraphics{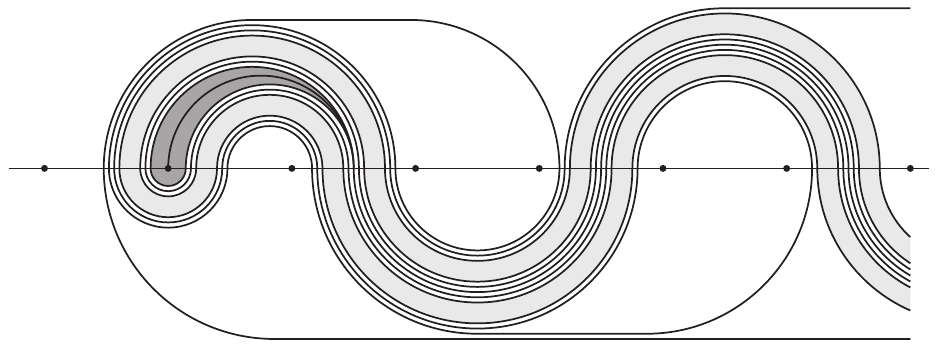}}
\put(0,47){$l$}
\end{picture}\\[5pt]
the leftmost part in the case where $n$ is even\\[12pt]
\begin{picture}(280,100)(0,0)
\put(7,0){\includegraphics{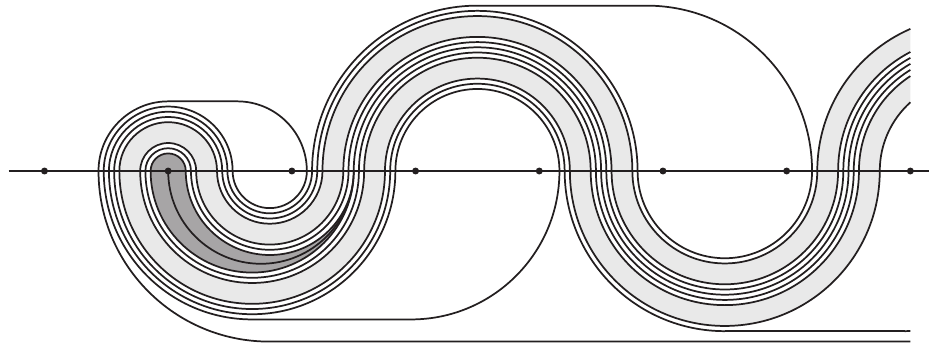}}
\put(0,47){$l$}
\end{picture}\\[5pt]
the leftmost part in the case where $n$ is odd\\[12pt]
\begin{picture}(280,100)(0,0)
\put(0,0){\includegraphics{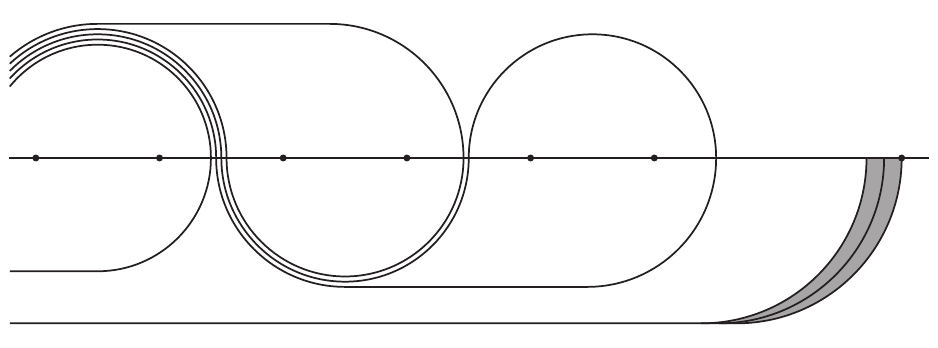}}
\put(272,51){$l$}
\end{picture}\\[5pt]
the rightmost part
\caption{The region $R(0)$ in the sphere $S_n$, whose long middle part looks like an arc.}
\label{fig_diagram_5}
\end{figure}

To draw the bridge diagram $(\sigma ^+,\sigma ^-)$ of $S_n$, we observe the deformation of the arcs of $\sigma ^+(s)$ as $s$ descends from $5$ to $0$.
Note that $\sigma ^+(5)$ is the collection of straight arcs in the sphere $S_n(5)$ as in Figure~\ref{fig_diagram_0}.
As $s$ descends from $5$ to $4$, the left endpoint of the rightmost arc $\sigma _n^+(s)$ moves as described by the gray arrow in the figure.
This gives us $\sigma ^+(4)\subset S_n(4)$ as in Figure~\ref{fig_diagram_1}.
By similar observation, as $s$ descends to $3$ (respectively, to $2$), we have $\sigma ^+(3)$ (respectively, $\sigma ^+(2)$) as in Figure~\ref{fig_diagram_2} (respectively, Figure~\ref{fig_diagram_3}).
Note that the region $R(2)\subset S_n(2)$ in Figure~\ref{fig_diagram_3} contains $n+1$ subarcs of $\sigma ^+(2)$ as in Figure~\ref{fig_tongue}.
We continue our observation as $s$ descends to $1$ and to $0$.
Though we do not include the whole picture of $\sigma ^+(1)$ (respectively, $\sigma ^+(0)$), by paying particular attention to the deformation of $R(2)$, we can find a region $R(1)\subset S_n(1)$ (respectively, $R(0)\subset S_n(0)$) as in Figure~\ref{fig_diagram_4} (respectively, Figure~\ref{fig_diagram_5}) which contains $n+1$ subarcs of $\sigma ^+(1)$ (respectively, $\sigma ^+(0)$) as in Figure~\ref{fig_tongue}.

\begin{figure}[p]
\begin{picture}(325,35)(0,0)
\put(5,5){\includegraphics{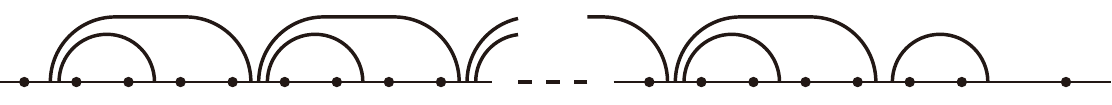}}
\put(0,7){$l$}
\put(32,0){$\delta _1$}
\put(62,0){$\delta _2$}
\put(92,0){$\delta _3$}
\put(122,0){$\delta _4$}
\put(208,0){$\delta _{n-3}$}
\put(238,0){$\delta _{n-2}$}
\put(268,0){$\delta _{n-1}$}
\put(317,0){$\delta _n$}
\end{picture}\\
in the case where $n$ is even\\[15pt]
\begin{picture}(325,35)(0,0)
\put(5,5){\includegraphics{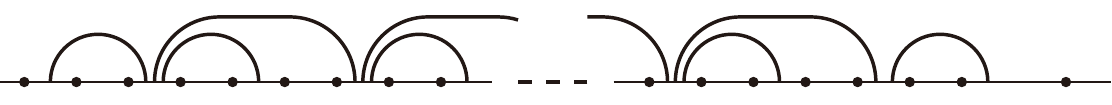}}
\put(0,7){$l$}
\put(32,0){$\delta _1$}
\put(62,0){$\delta _2$}
\put(92,0){$\delta _3$}
\put(122,0){$\delta _4$}
\put(208,0){$\delta _{n-3}$}
\put(238,0){$\delta _{n-2}$}
\put(268,0){$\delta _{n-1}$}
\put(317,0){$\delta _n$}
\end{picture}\\
in the case where $n$ is odd
\caption{Some components of $R(0)\cap S_+$.}
\label{fig_diagram_6}
\end{figure}

We can see that the pair $\left( (\sigma ^+,\sigma ^-),l\right) $ satisfies the $2$-connected condition as follows.
We choose $S_+$ (respectively, $S_-$) to be the hemisphere above (respectively, below) $l$ in Figure~\ref{fig_diagram_5}, and let $\delta _1,\delta _2,\ldots ,\delta _n$ and ${\mathcal G}_{i,j,\varepsilon }$ be as in Section~\ref{diagrams}.
We can see from Figure~\ref{fig_diagram_5} that there exist components of $R(0)\cap S_+$ as in Figure~\ref{fig_diagram_6}.
In particular, there exists at least one component of $R(0)\cap S_+$ separating $\delta _i$ from $\delta _j$ in $S_+$ for each pair of distinct $i,j\in \{ 1,2,\ldots ,n\} $.
Since $R(0)$ contains subarcs of $\sigma ^+(0)$ as in Figure~\ref{fig_tongue}, the component of $R(0)\cap S_+$ contains a subarc of $\sigma _v^+$ and a subarc of $\sigma _{v+1}^+$ that are adjacent in $S_+$ for each $v\in \{ 1,2,\ldots ,n\} $, where the index $v+1$ is considered modulo $n$.
The graph ${\mathcal G}_{i,j,+}$ therefore has a cycle passing through the vertices $1,2,\dots,n$ in this order, which implies that ${\mathcal G}_{i,j,+}$ is $2$-connected.
Similarly, one can also check that ${\mathcal G}_{i,j,-}$ is $2$-connected for each pair of distinct $i,j\in \{ 1,2,\ldots ,n\} $.

By Theorem~\ref{criterion}, the $n$-bridge sphere $S_n$ is strongly irreducible, and hence destabilized.
It follows that the knot $K_n$ is not a rational knot, since a rational knot cannot admit a destabilized $n$-bridge sphere by the results mentioned in Section~\ref{introduction}.
It follows that the $3$-bridge sphere $S_n^\perp $ is also destabilized, since only rational knots may admit a $2$-bridge sphere.
This completes the proof of Theorem~\ref{biplat}.

\end{document}